\documentclass[10pt]{article}
\usepackage{amsmath}
\usepackage{latexsym}
\usepackage{amssymb}
\usepackage{amsfonts}
\usepackage{amsthm}
\title{ A MATHEMATICAL COMMENT ON  \\  LANCZOS POTENTIAL THEORY }
\author{J.-F. Pommaret \\ CERMICS, Ecole des Ponts ParisTech, France \\
jean-francois.pommaret@wanadoo.fr \\
( http://cermics.enpc.fr/$\sim$pommaret/home.html) }
\date{  }
\begin{document}
\maketitle

\noindent
{\bf ABSTRACT}  \\

The last invited lecture published in $1962$ by Lanczos on his potential theory is never quoted because it is in french. Comparing it with 
a commutative diagram in a recently published paper on gravitational waves, we suddenly understood the confusion made by Lanczos between Hodge duality and differential duality. Our purpose is thus to revisit the mathematical framework of Lanczos potential theory in the light of this comment, getting closer to the formal theory of Lie pseudogroups through differential double duality and the construction of finite length differential sequences for Lie operators. We use the fact that a differential module $M$ defined by an operator ${\cal{D}}$ with coefficients in a differential field $K$ has vanishing first and second differential extension modules if and only if its adjoint differential module $N=ad(M)$ defined by the adjoint operator $ad({\cal{D}})$ is reflexive, that is $ad({\cal{D}})$ can be parametrized by the operator $ad({\cal{D}}_1)$ when ${\cal{D}}_1$ generates the compatibilty conditions (CC)
 of ${\cal{D}}$ while $ad({\cal{D}}_1)$ can be parametrized by $ad({\cal{D}}_2)$ when ${\cal{D}}_2$ generates the CC of ${\cal{D}}_1$. We provide an explicit description of the potentials allowing to parametrize the Riemann and the Weyl operators in arbitrary dimension, both with their respective adjoint operators.   \\

\vspace{3cm}


\newpage

\noindent
{\bf 1) INTRODUCTION}\\

We start recalling that a short exact sequence of modules  $0 \rightarrow M' \stackrel{f}{\longrightarrow} M \stackrel{g}{\longrightarrow} M" \rightarrow 0 $  {\it splits} if and only if $M\simeq M' \oplus M"$ or, equivalently, if and only if there exists an epimorphism $u:M \rightarrow M'$ such that $u \circ f=id_{M'}$ or a monomorphism $v:M" \rightarrow M$ such that $g\circ v \simeq id_{M"}$ ([35], p 33). In Riemannian geometry, we have the classical formula relating the (linearized) Riemann $4$-tensor $R$ and the (linearized) Weyl $4$-tensor $W$ in arbitrary dimension $n\geq 3$ when $\omega$ is a non-degenerate metric with constant Riemannian curvature or flat like the Minkowski metric ([33]):  \\
\[   W_{kl,ij}=R_{kl,ij}  - \frac{1}{(n-2)}({\omega}_{ki}R_{lj}-{\omega}_{kj}R_{li} + {\omega}_{lj}R_{ki} - {\omega}_{li}R_{kj}) +
 \frac{1}{(n-1)(n-2)}({\omega}_{ki}{\omega}_{lj} - {\omega}_{kj}{\omega}_{li})R   \] 
We have proved in ([28],[32]) the existence of an {\it intrinsic splitting} of the short exact sequence of tensor bundles  $0 \rightarrow Ricci \stackrel{f}{\longrightarrow} Riemann \stackrel{g}{\longrightarrow} Weyl  \rightarrow 0$ where $f$ is defined by equaling the right member to zero while defining $u$ to be $(R^k_{l,ij}) \rightarrow (R^r_{i,rj}=R_{ij}=R_{ji})$ in order to obtain an isomorphism $Riemann \simeq Ricci \oplus Weyl$ that {\it cannot} be understood without the use of the Spencer $\delta$-cohomology. This result, that will be revisited in Sections $3$ and $4$, is coherent with the formulas $(15)$ and $(16)$ of ([13]) when $n=4$ if we notice that $U+V=2R$ and set $V=2W$ .\\

Introducing the {\it Lanczos potential} $L=(L_{ij,k})$ as a $3$-tensor satisfying the algebraic relations:  \\
\[  \fbox{ $L_{ij,k} + L_{ji,k}=0, \,\,\,\,\,  L_{ij,k} + L_{jk,i} + L_{ki,j} = 0 $} \]
Lanczos claimed in the formula $(III.5)$ of ([11]) or $(17)$ of ([13]) to have parametrized the Riemann tensor through the relation:  \\
\[   R_{kl,ij} = {\nabla}_j L_{kl,i} - {\nabla}_i L_{kl,j} + {\nabla}_l L_{ij,k} - {\nabla}_k L_{ij,l}  \]
where $\nabla$ is the covariant derivative. However, even if we can easily verify the algebraic conditions that {\it must} be satisfied by a {\it Riemann candidate} with $4$ indices, namely:  \\
\[    R_{kl,ij}= - R_{lk,ij} = - R_{kl,ji} = R_{ij,kl}, \,\,\,\,   R_{kl,ij} + R_{ki,jl} + R_{kj,li} = 0   \]
the generating {\it compatibility conditions} (CC) of the underlying operator for the left member {\it cannot} be the {\it Bianchi identities}:  \\
\[  {\nabla}_r{R^k}_{l,ij} + {\nabla}_i{R^k}_{l,jr} + {\nabla}_j{R^k}_{l,ri}= 0  \]
which are produced by the well known parametrization described by the $Riemann$ operator acting on a perturbation $\Omega \in S_2T^*$ of the background metric $\omega$, that is, when $\omega$ is the Minkowski metric:  \\
\[  2R_{kl,ij}=(d_{li}{\Omega}_{kj} - d_{lj}{\Omega}_{ki}) - (d_{ki}{\Omega}_{lj} - d_{kj}{\Omega}_{li}) \,\,\, \Rightarrow  \,\,\,  
d_r{R^k}_{l,ij} + d_i{R^k}_{l,jr} + d_j{R^k}_{l,ri}= 0  \]  
This contradiction can also been checked directy by substitution because we have:  \\
\[ R_{kl,ij} = d_j L_{kl,i} - d_i L_{kl,j} + d_l L_{ij,k} - d_k L_{ij,l} \,\,\, \Rightarrow \,\,\, d_r{R^k}_{l,ij} + d_i{R^k}_{l,jr} + d_j{R^k}_{l,ri}
\neq 0   \]
Setting now $A_{ik}=A_{ki}= \frac{1}{2} ({\nabla}_r {{L_i}^r}_{,k }+{\nabla}_r{{L_k}^r}_{,i})$, we obtain the so-called Weyl-Lanczos equations given in ([13], formula (17)):   \\
\[  \begin{array}{rcl}
W_{kl,ij} &  =  &  {\nabla}_j L_{kl,i}  - {\nabla}_i L_{kl,j} + {\nabla}_l L_{ij,k} - {\nabla}_k L_{ij,l}  \\
              &   &   + {\omega}_{li}A_{jk} - {\omega}_{jl}A_{ik} + {\omega}_{jk}A_{il} - {\omega}_{ik}A_{jl}                                                
\end{array}   \]    \\
and it becomes clear that Lanczos was not even aware of the Weyl tensor at that time ([15]), only knowing the algebraic conditions that must be fulfilled by $C$, namely:  \\   
\[    W_{kl,ij}= - W_{lk,ij} = - W_{kl,ji} = W_{ij,kl}, \,\,\,\,   W_{kl,ij} + W_{ki,jl} + W_{kj,li} = 0 , \,\,\, {W^r}_{i,rj}= 0  \]
the last condition reducing the number of linearly independent components from $20$ to $10$ for space-time, that is when the dimension is $n=4$. As a byproduct, the previous contradiction still holds. \\
In order to recapitulate the above procedure while setting $F_0=S_2T^*$, we have the differential sequence:  \\
\[  0 \rightarrow \Theta \rightarrow T \stackrel{Killing}{\longrightarrow} F_0 \stackrel{Riemann}{\longrightarrow} F_1 \stackrel{Bianchi}{\longrightarrow} F_2 \rightarrow ...  \]
\[0 \rightarrow \Theta \rightarrow n \rightarrow n(n+1)/2 \rightarrow n^2(n^2-1)/12 \rightarrow n^2(n^2-1)(n-2)/24 \rightarrow ...  \]
where $\Theta$ is the sheaf of Killing vector fields for the Minkowski metric. \\
For historical reasons, defining the operators $Cauchy=ad(Killing)$, $Beltrami=ad(Riemann)$ and  $Lanczos=ad(Bianchi)$, we obtain the adjoint sequence:  \\
\[  0 \leftarrow ad(T) \stackrel{Cauchy}{\longleftarrow} ad(F_0) \stackrel {Beltrami}{\longleftarrow} ad(F_1)\stackrel{Lanczos}{\longleftarrow} ad(F_2) \leftarrow ...   \]
where $ad(E)={\wedge}^nT^*\otimes E^*$ for any vector bundle $E$ with $E^*$ obtained from $E$ by inverting the transition rules when changing local coordinates, exactly like $T \rightarrow (T)^*=T^*$. Accordingly, all the problem will be to prove that {\it each operator is indeed parametrized by the preceding one}. As we shall see, the conformal situation could be treated similarly while starting with the {\it conformal Killing} operator followed by the $Weyl$ operator and replacing each {\it classical} vector bundle $F$ by the corresponding {\it conformal} bundle $\hat{F}$. However, such a point of view is leading to a {\it true nonsense} because we shall discover that the analogue of the Bianchi operator is of order $2$, ...  just when $n=4$. This striking result has been confirmed by computer algebra and the reader can even find the details in book form ([29]). It follows that both the Riemann and Weyl frameworks of the Lanczos potential theory must be entirely revisited.   \\

 C. Lanczos (1893-1974) wrote a book on variational calculus ([11]) and three main papers (1939,1949,1962) on the {\it potential theory} in physics, mostly by comparing the case of electromagnetism (Maxwell equations) with the search for parametrizing the Riemann and Weyl tensors ([9],[10],[12]) and we refer the reader to the nice historical survey ([15]) for more details. However, Lanczos has been invited in 1962 by Prof. A. Lichnerowicz to lecture in France and this lecture has been published in french ([13]). He got inspiration from what happens in electromagnetism (EM) where the first set of Maxwell equations $dF=0$ is saying that the EM field $F\in {\wedge}^2T^*$ is a closed $2$-form that can be parametrized by $dA=F$ for an arbitrary EM potential $A\in T^*$  ([31]). Accordingly, Lanczos created the concept of "{\it candidate} " while noticing that the Riemann and Weyl $4$-tensors must "{\it a priori} " satisfy algebraic relations reducing the number of their components $R_{kl,ij}$ and $W_{kl,ij}$ respectively to $20$ and $10$ when $n=4$ as we saw. Now, we have proved in many books ([16-19],[21],[29],[32]) or papers that it is not possible to understand the mathematical structure of the Riemann and Weyl tensors, both with their splitting link, without the following four important comments:  \\

\noindent
$\bullet$ The clever results discovered by E. Vessiot as early as in $1903$ ([38]) are still neither known nor acknowledged today, though they generalize the constant Riemaniann curvature integrability condition discovered $25$ years later by L. P. Eisenhart ([7]). They also allow to understand the direct link existing {\it separately} between the Riemann tensor and the Lie group of isometries (considered as a Lie pseudogroup) of a non-degenerate metric on one side or between the Weyl tensor and the Lie group of conformal isometries (considered again as a Lie pseudogroup) of this metric on another side. With more details, Vessiot proved that, for any Lie pseudogroup $\Gamma \in aut(X)$ one can find a {\it geometric object } $\omega$, may be of a high order $q$ and not of a tensorial nature, which is characterizing $\Gamma$ as the group of local diffeomorphisms preserving $\omega$, namely:  \\
\[\    \Gamma = \{f\in aut(X)\mid {\Phi}_{\omega}(j_q(f))=j_q(f)^{-1}(\omega)=\omega \}   \]   
where $\omega$ must satisfy certain (non-linear in general) integrability conditions of the form:   \\
\[ I(j_1(\omega))=c(\omega)\]
when $q$ is large enough, called {\it Vessiot structure equations}, locally depending on a certain number of constants, now called {\it Vessiot structure constants}, and we let the reader compare this situation to the Riemann case ([25],[29],[32]). These structure equations were perfectly known by E. Cartan (1869-1951) who {\it never} said that they were at least competing with or even superseding the corresponding {\it Cartan structure equations } that he developped about at the same time for similar purposes. The underlying reason is of a purely personal origin related to the origin of  "{\it differential Galois Theory} " within a kind of "{\it mathematical affair} " involving the best french mathematicians of that time (H. Poincar\'{e}, E. Picard, G. Darboux, P. Painlev\'{e}, E. Borel, ...). The main original letters, given to the author of this paper by M. Janet, have ben published in ([18]) and can now be examined in the main library of Ecole Normale Sup\'{e}rieure in Paris where they have been deposited. \\

\noindent
$\bullet$ A nonlinear operator with second member does not in general admit CC, ... unless it corresponds to the defining equations in Lie form of a Lie pseudogroup and the CC are the Vessiot structure equations in that case with structure constants determined by the chosen geometric object (compare again to the Riemannian geometry). We have shown in many books already quoted that, if ${\cal{D}}$ is a Lie operator, that is $[\Theta, \Theta]\in \Theta $ with bracket induced by the ordinary bracket of vector fields, then the system ${\cal{D}}\xi=\Omega$ is the linearization of a non-linear version when $\Omega$ is a perturbation of $\omega$ (twice the infinitesimal deformation tensor in elasticity) along the formula:  \\
\[  {\cal{D}}\xi= {\cal{L}}(\xi) \omega= \frac{d}{dt} (j_q(exp(t\xi))^{-1} (\omega)){\mid}_{t=0}  \]             
Similarly, we can choose for the corresponding generating CC ${\cal{D}}_1$ the linearization of a non-linear version described by the Vessiot structure equations:  \\
\[     \frac{\partial I}{\partial j_1(\omega)}(j_1(\omega)) j_1(\Omega) =\frac{\partial c}{\partial \omega}(\omega)\Omega   \]
that is {\it exactly} what we did for the flat Minkowski metric. However, Lanczos has been studying the CC ${\cal{D}}_2$ of ${\cal{D}}_1$, ignoring that, contrary to the previous situation, ${\cal{D}}_2$ almost never comes from a linearization. It is therefore quite strange to discover that Lanczos {\it never} discovered that what he was doing with ${\cal{D}}_1$ and ${\cal{D}}_2$ while using quadratic Lagrangians in $R$ along ([9-13]), was {\it exactly} what is done in any textbook of elasticity or continum mechanics with ${\cal{D}}$ and ${\cal{D}}_1$ while using quadratic Lagrangians in 
$\Omega$ ([29],[30],[32]). We believe that Lanczos was too much obsessed by comparing $R$ in GR to $F$ in EM and we shall provide more details at the end of Section 2.        \\

\noindent
$\bullet$ We now present the main origin of the troubles met by Lanczos and followers. With only a few words, we may say that most physicists and even many mathematicians are not familiar with the modern developments of differential homological algebra, in the sense that they do believe, starting from a linear differential operator ${\cal{D}}: E \rightarrow F=F_0:\xi \rightarrow \eta$ between (the sections of) two given vector bundles, one may construct ({\it at least}) the generating CC as an operator ${\cal{D}}_1:F_0 \rightarrow F_1: \eta \rightarrow \zeta$ and so on. As a byproduct, knowing only this " {\it step by step} " procedure, they are largely unaware of a " {\it global} " procedure, apart from the very specific case of the Poincar\'{e} sequence for the exterior derivative $d: {\wedge}^rT^* \rightarrow {\wedge}^{r+1}T^*$. Accordingly, $E$ and $F_0$ being given, they are loosing any relationship that could exist with $F_1,F_2, ..., F_n$. For example, the author of this paper perfectly remembers that, when he was a student of Prof. A. Lichnerowicz in the $70^{es}$, he was hardly able to know about the number $n^2(n^2-1)/12$ of components of the Riemann tensor found by E. Cartan with tedious combinatorics but totally unable to know the number $n^2(n^2-1)(n-2)/24$ of Bianchi identities. We find therefore useful to recall a few historical facts because this last number, namely $20=24 - 4$, is {\it exactly} the number of {\it Lagrange multipliers} used by Lanczos in his variational approach to Riemanniann geometry ([13]).  \\

The first {\it finite length differential sequence} (now called {\it Janet sequence}) has been exhibited as a footnote by M. Janet in 1920 ([8],[16]) with only $E,F_0, F_1, ... ,F_n$ when ${\cal{D}}$ is involutive and $n$ is the number of independent variables. It must be noticed that the {\it exterior calculus} of E. Cartan ([1]), mixing up the dependent {\it and} independent variables, has put a shadow on this point of view, also combined with the non-intrinsic approach for finding Gr\"{o}bner bases along similar methods. It is only during the period $1960-1970$ that D.C. Spencer and coworkers brought new differential homological algebraic intrinsic methods for studying such sequences ([37]). Then analysis became such a fashionable subject that almost nobody took the risk to use these difficult new methods in physics, {\it even though they are largely superseding the previous ones}. In order to convince the reader about the problems that could be met, we end this comment with two unusual examples.  \\

\noindent
{\bf EXAMPLE 1.1}: Using standard notations from jet theory, let us consider the trivially involutive operator $j_q: E \rightarrow J_q(E):({\xi}^k) \rightarrow ({\xi}^k, {\partial}_i{\xi}^k, {\partial}_{ij}{\xi}^k, ... )$ up to order $q$. The symbol $g_q$ is trivially involutive because it is defined by the linear system ${\xi}^k_{i_1 ... i_q}= 0$ and thus $g_q=0$. On one side, setting $J_q(E)=F$ and $F=F_0$, the (canonical) Janet sequence {\it must} be of the form:  \\
\[  0 \rightarrow E \stackrel{{\cal{D}}}{\longrightarrow} F_0 \stackrel{{\cal{D}}_1}{\longrightarrow} ... \stackrel{{\cal{D}}_n}{\longrightarrow}
 F_n \rightarrow 0  \]
in which ${\cal{D}}=j_q$. However, on the other side, setting $J_q(E)=C_0(E)$ and exhibiting the {\it Spencer} bundles $C_r(E)={\wedge}^rT^*\otimes J_q(E)/ \delta ({\wedge}^{r-1}\otimes S_{q+1}T^*\otimes E)$ with the {\it Spencer} operator ${\cal{D}}_{r+1}: C_r(E) \rightarrow C_{r+1} (E)$ induced by the (naive) Spencer operator $D$ and its extensions:  \\
\[   D: J_{q+1}(E) \rightarrow T^*\otimes J_q(E): {\xi}_{q+1} \rightarrow j_1({\xi}_q) -{\xi}_{q+1}: {\xi}_{q+1} \rightarrow ({\partial}_i{\xi}^k-{\xi}^k_i, {\partial}_i{\xi}^k_j - {\xi}^k_{ij}, ...)\]
\[  D: {\wedge}^rT^*\otimes J_{q+1}(E) \rightarrow {\wedge}^{r+1}T^* \otimes J_q(E): \alpha \otimes {\xi}_{q+1}\rightarrow  d\alpha \otimes {\xi}_q +{(-1)}^r \alpha \wedge D{\xi}_{q+1}  \]
We have $D\circ D=0\Rightarrow D_{r+1}\circ D_r=0$ and, because of the exterior product, we should obtain the finite length {\it second Spencer sequence}:  \\
\[  0 \rightarrow E \stackrel{j_q}{\longrightarrow} C_0(E) \stackrel{D_1}{\longrightarrow} C_1(E) \stackrel{D_2}{\longrightarrow} ... \stackrel{D_n}{\longrightarrow}C_n(E) \rightarrow 0   \]
We point out the fact that, in the sequence:  \\
\[     J_{q+1}(E)\stackrel{D}{\longrightarrow} T^*\otimes J_q(E) \stackrel{D}{\longrightarrow}{\wedge}^2T^*\otimes J_{q-1}(E)  \]
the second operator is only generating the CC of the first but not describing them totally because ${\partial}_j({\partial}_i{\xi}^k- {\xi}^k_i) - {\partial}_i({\partial}_j{\xi}^k-{\xi}^k_j)= {\partial}_i{\xi}^k_j - {\partial}_j{\xi}^k_i=0$. \\
Comparing the two methods, in the Janet sequence we have a "{\it step by step} " procedure while, in the Spencer sequence we have a " {\it global} " procedure. Before looking at the next lines, we invite the reader to understand how difficult is the computation of $C_r(E)$ and $F_r$ even when $n=3,m=1,q=2$ and to conclude that ... $C_r(E)\simeq F_r$ (!). \\
Let us consider the general situation of an involutive operator ${\cal{D}}=\Phi \circ j_q$ defined by an involutive system $R_q=ker(\Phi)$. As in ([16],[19],[21]), we may define $F_0$ by the short exact sequence $0 \rightarrow R_q \rightarrow J_q(E) \stackrel{\Phi}{\longrightarrow} F_0 \rightarrow 0  $ and the {\it Janet bundles} $F_r$ by the short exact sequences:  \\
\[  0 \rightarrow {\wedge}^rT^*\otimes R_q + \delta({\wedge}^{r-1}T^*\otimes S_{q+1}T^*\otimes E ) \rightarrow  {\wedge}^rT^*\otimes J_q(E) \rightarrow  F_r \rightarrow 0   \]
where "$+$" denotes a sum that may not be a direct sum. We may define similarly the {\it Spencer bundles} $C_r$ by the short exact sequences: \\
\[    0 \rightarrow \delta({\wedge}^{r-1}T^* \otimes g_{q+1}) \rightarrow {\wedge} ^rT^*\otimes R_q \rightarrow C_r  \rightarrow 0  \]
and we have the short exact sequences $0 \rightarrow C_r \rightarrow C_r(E) \rightarrow F_r \rightarrow 0$ showing that, in general, {\it the Janet and Spencer sequences are two completely different differential sequences} used as resolutions of the sheaf $\Theta$ of solutions of ${\cal{D}}$. However, {\it in the present situation}, we have by definition $F_0=J_q(E) \Rightarrow R_q=0 \Rightarrow F_r=C_r(E)$ as an unexpected result involving the Spencer operator. For a later use, we notice that any inclusion $R_q \subset {\hat{R}}_q$ of involutive systems provides canonical monomorphisms $0 \rightarrow C_r\rightarrow {\hat{C}}_r$ both with canonical epimorphisms $F_r \rightarrow {\hat{F}}_r \rightarrow 0$, a result showing that the Janet and Spencer sequences cannot be used equivalently in physics, even if they can be used equivalently in mathematics for computing differential extension modules (See [23, [25] and [26] for the best examples we know concerning elasticity, general relativity and gauge theory).  \\

\noindent
{\bf EXAMPLE 1.2}: In Riemannian geometry, the situation is even more tricky, a {\it fact} explaining why all the tentatives made by the various authors ([2-6]) in order to use either Cartan ([3]) or Janet ([2]) or Gr\"obner are {\it never} appealing to differential sequences. They could not therefore describe any link between the work of Lanczos and the construction of differential sequences, in particular the close relationship existing between Lanczos Lagrange multipliers and Spencer $\delta$-cohomology that we shall exhibit in Sections $2$ and $3$. Moreover, as any action of a Lie group has a finite number of infinitesimal generators $\{ {\theta}_{\tau} \}$ providing a Lie algebra ${\cal{G}}$, we have $C_r={\wedge}^rT^*\otimes R_q\simeq {\wedge}^rT^*\otimes {\cal{G}}$ when $q$ is large enough. In the particular case of isometries and conformal isometries, we have $q=1$ or $q=2$ and this number is respectively equal to $n(n+1)/2$ and $(n+1)(n+2)/2$. However, the corresponding canonical Janet and Spencer sequences cannot be constructed for $R_1$ or ${\hat{R}}_1$ which are {\it not} involutive but {\it must} be constructed for $R_2$ or ${\hat{R}}_3$ ({\it care}) which are involutive ([33]). It will follow that the Janet and Spencer sequences are completely different and {\it must} be therefore carefully distinguished as we shall discover that the step by step construction of generating CC will bring a lot of surprises for the successive operators ${\cal{D}}, {\cal{D}}_1,{\cal{D}}_2$ or $\hat{\cal{D}}, {\hat{\cal{D}}}_1, {\hat{\cal{D}}}_2$.  \\ \\

\noindent
$\bullet$ {\it Last but not least}, we explain the way towards the solution of the parametrization problem by means of {\it differential double duality}. 
The starting point has been a challenge proposed in 1970 by J. Wheeler to find out a parametrization of Einstein equations in vacuum, at the time the 
author of this paper was a student of D.C. Spencer in Princeton university. Later on, discovering by chance while teaching elasticity, that the parametrization of Cauchy equations in dimension $2$ by the Airy function was nothing else than the formal adjoint of the {\it Riemann} operator, this result allowed him to give a {\it negative answer} to the challenge in 1995 and we point out that not a single step ahead had been produced during the previous $25$ years ([20]). Then, again by chance, we discovered in 1995 the english translation of a thesis by M. Kashiwara  (See [36] for a more accessible version). This has been the starting point of the use of {\it differential modules} and {\it differential homological algebra} for applications ([20],[32]), thanks to the pioneering work of U.Oberst in 1990 on control theory where a control system is controllable if and only if it is parametrizable ([14],[21],[22]). Finally, the author thanks Prof. L. Andersson met in Paris during the Einstein centenary, still by chance, to have suggested him to apply these new techniques in order to study the Lanczos potential theory. \\

\vspace{2cm}

\noindent
{\bf 2) RIEMANN/LANCZOS POTENTIAL}  \\

Having in mind the variational procedure used in optimal control theory when $n=1$ and in EM when $n=4$, let us assume that the differential sequence:  \\
\[ \xi  \stackrel{\cal{D}}{\longrightarrow} \eta  \stackrel{{\cal{D}}_1}{\longrightarrow}  \zeta   \]
is {\it formally exact}, that is ${\cal{D}}_1$ generates the CC of ${\cal{D}}$ and thus $\xi$ is a potential for ${\cal{D}}_1$. We may consider a variational problem for a cost function or lagrangian $\varphi (\eta)$ under the linear OD or PD constraint described by ${\cal{D}}_1\eta=0$. \\

\noindent
$\bullet$ Introducing convenient Lagrange multipliers $\lambda$ while setting $dx=dx^1\wedge ... \wedge dx^n$ for simplicity, we must vary the integral:  \\
 \[       \Phi=\int [\varphi(\eta) - \lambda {\cal{D}}_1\eta]dx  \Rightarrow \delta \Phi=\int [(\partial\varphi(\eta)/\partial\eta)\delta\eta - \lambda{\cal{D}}_1\delta\eta]dx \]
Integrating by parts, we obtain the {\it Euler-Lagrange} (EL) equations: \\
\[  \partial\varphi(\eta)/\partial\eta = ad({\cal{D}}_1)\lambda  \]
to which we have to add the constraint ${\cal{D}}_1\eta=0$ obtained by varying $\lambda$ independently. If $ad({\cal{D}}_1)$ is an injective operator, in particular if ${\cal{D}}_1$ is formally surjective (no CC) while $n=1$ as in OD optimal control and the differential module defined by ${\cal{D}}_1$ is torsion-free, thus free ([13],[28]) or $n\geq 1$ and this module is projective, then {\it one can obtain} $\lambda$ {\it explicitly} and eliminate it by substitution. Otherwise, using the generating CC ${\cal{D}}'$ of $ad({\cal{D}}_1)$, we have to study the formal integrability of the combined system: \\
\[    {\cal{D}}'\partial\varphi(\eta)/\partial\eta=0, \hspace{4mm} {\cal{D}}_1\eta=0   \]
which may be a difficult task as can be seen in the examples of the Introduction of ([28],[34]).   \\

\noindent
$\bullet$ However, {\it we may also transform the given variational problem with constraint into a variational problem without any constraint if and 
only if the differential constraint can be parametrized}. \\
Indeed, using the parametrization of ${\cal{D}}_1$ by ${\cal{D}}$, we may vary the integral: \\
\[  \Phi=\int \varphi({\cal{D}}\xi)dx  \Rightarrow \delta \Phi = \int (\partial\varphi(\eta)/\partial\eta){\cal{D}}\delta\xi  dx                \]
whenever $\eta={\cal{D}}\xi$ and integrate by parts for arbitrary $\delta\xi$ in order to obtain the EL equations:  \\
\[   ad({\cal{D}})\partial\varphi(\eta)/\partial\eta=0, \hspace{4mm}  \eta={\cal{D}}\xi  \] 
 in a way which is coherent with the previous approach {\it if and only if} $ad({\cal{D}})$ generates the CC of $ad({\cal{D}}_1)$, a condition rarely satisfied in general. \\
 
Accordingly, even if ${\cal{D}}_1$ generates the CC of ${\cal{D}}$, in general $ad({\cal{D}})$ may not generate all the CC of $ad({\cal{D}}_1)$ in the adjoint differential sequence:  \\
\[    \nu \stackrel{ad(\cal{D})}{\longleftarrow}  \mu  \stackrel{ad({\cal{D}}_1)}{\longleftarrow}  \lambda    \]
Such a striking "{\it gap} ", namely the lack of formal exactness of the adjoint sequence when the initial sequence is formally exact, led to introduce the {\it differential extension modules} because of the following (difficult) theorems (See [21,35] or [25,28,30,32] for more details):  \\

\noindent
{\bf THEOREM 2.1}: If $M$ is the differential module defined by ${\cal{D}}$, the {\it extension modules} $ext^i(M)$ do not depend on the sequence used for their computation and are {\it torsion modules} for $i\geq 1$, that is to say $rk_D(ext^i(M))=0, \forall i\geq 1$.  \\

With the same notations, let us introduce the differential module $N=ad(M)$ defined by the adjoint operator $ad({\cal{D}})$. \\

\noindent
{\bf THEOREM 2.2}: When $M={ }_DM$ is a (left) differential module, then $hom_D(M,D)$ is a right differential module because $D$ is a bimodule over itself and $N={ }_DN=hom_K({\wedge}^nT^*, hom_{D}(M,D))$ is again a (left) differential module.  \\

\noindent
{\bf COROLLARY 2.3}: A differential module is {\it torsion-free} if and only if it can be embedded into a free differential module.   \\

\noindent
{\bf COROLLARY 2.4}: The differential module $M$ is such that $ext^1(M)=0, ext^2(M)=0$ if and only if the differential module $N$ is {\it reflexive}, that is $ad({\cal{D}})$ can be parametrized by $ad({\cal{D}}_1)$ which can be itself parametrized by $ad({\cal{D}}_2)$. When $n=3$, the simplest example is the {\it div} operator that can be parametrized by the {\it curl} operator which can be parametrized by the {\it grad} operator.  \\
 
\noindent
{\bf THEOREM 2.5}: The {\it Spencer sequence} for any Lie operator ${\cal{D}}$ which is coming from a Lie group of transformations, with a Lie group $G$ acting on $X$, is (locally) isomorphic to the tensor product of the Poincar\'{e} sequence by the Lie algebra ${\cal{G}}$ of $G$.  \\

\noindent
{\it Proof}: If $M$ is the differential module defined by ${\cal{D}}$, we want to prove that the extension modules $ext^1(M)$ and $ext^2(M)$ vanish, that is, if ${\cal{D}}_1$ generates the CC of ${\cal{D}}$ but also ${\cal{D}}_2$ generates the CC of ${\cal{D}}_1$, then $ad({\cal{D}})$ generates the CC of $ad({\cal{D}}_1)$ and $ad({\cal{D}}_1)$ generates the CC of $ad({\cal{D}}_2)$. We also remind the reader that we have shown in ([28],[34]) that it is not easy to exhibit the CC of the Maxwell or Morera parametrizations when $n=3$ and that a direct checking for $n=4$ should be strictly impossible. It has been proved by L. P. Eisenhart in 1926 ([7]) that the solution space $\Theta$ of the Killing system has $n(n+1)/2$ {\it infinitesimal generators} $\{ {\theta}_{\tau}\}$ linearly independent over the constants if and only if $\omega $ had constant Riemannian curvature, namely zero in our case. As we have a transitive Lie group of transformations preserving the metric considered as a transitive Lie pseudogroup, the three classical theorems of Sophus Lie assert than $[{\theta}_{\rho},{\theta}_{\sigma}]=c^{\tau}_{\rho \sigma} {\theta}_{\tau}$ where the {\it structure constants} $c$ define a Lie algebra ${\cal{G}}$. We have therefore $\xi \in \Theta \Leftrightarrow \xi = \{{\lambda}^{\tau}{\theta}_{\tau}\mid {\lambda}^{\tau}=cst\}$. Hence, we may replace locally the Killing system by the system ${\partial}_i{\lambda}^{\tau}(x)=0$, getting therefore the differential sequence:  \\
\[  0 \rightarrow \Theta \rightarrow {\wedge}^0T^*\otimes {\cal{G}} \stackrel{d}{\longrightarrow} {\wedge}^1T^*\otimes {\cal{G}} \stackrel{d}{\longrightarrow} ... \stackrel{d}{\longrightarrow} {\wedge}^nT^* \otimes {\cal{G}} \rightarrow 0  \]
 which is the tensor product of the Poincar\'{e} sequence by ${\cal{G}}$. Finally, it follows from the above Theorem that the extension modules considered do not depend on the resolution used and thus vanish because the Poincar\'{e} sequence is self adjoint (up to sign), that is $ad(d)$ generates the CC of $ad(d)$ at any position, exactly like $d$ generates the CC of $d$ at any position. This (difficult) result explains why the adjoint differential modules we shall meet will be torsion-free or even reflexive. We invite the reader to compare with the situation of the Maxwell equations in electromagnetisme (See [18, p 492-494] for more details). However, we have explained in ([17],[24],[26],[27],[32]) why neither the Janet sequence nor the Poincar\'{e} sequence can be used in physics and must be replaced by the {\it Spencer sequence} which is another resolution of $\Theta$ ([19],[32],[37]). Though this is out of the scope of this paper, we shall nevertheless shortly describe the relation existing between the above results and the Spencer operator, thus the Spencer sequence. For this, let us define for any $q\geq 0$ the section ${\xi}_q={\lambda}^{\tau}(x)j_q({\theta}_{\tau})(x)=({\xi}^k_{\mu}(x)={\lambda}^{\tau}(x) {\partial}_{\mu}{\theta}^k_{\tau}(x))\in R_q$. With the standard notations of ([16],[19],[21]) and $0 \leq \mid \mu \mid \leq q$, the components of the Spencer operator become:  \\
 \[   D{\xi}_{q+1}=j_1({\xi}_q )- {\xi}_{q+1}= ({\partial}_i{\xi}^k_{\mu}(x)- {\xi}^k_{\mu + 1_i}(x)) =(({\partial}_i{\lambda}^{\tau}(x)){\partial}_{\mu}{\theta}^k_{\tau}(x) {\xi})\in T^*\otimes R_q\]
When $q$ is large enough, that is $q=1$ for the Killing and conformal Killing systems, we obtain the desired identification justifying our claim.  \\
\hspace*{12cm}   Q.E.D.  \\

\noindent
{\bf COROLLARY 2.6}: When ${\cal{D}}$ is the Killing operator, then $ {ext}^1(M)=0,\,\, {ext}^2(M)=0 $ and there is no gap. The situation is similar if we start with the conformal Killing operator $\hat{{\cal{D}}}$.  \\

\noindent
{\bf REMARK 2.7}: If the differential module $M$ defined by ${\cal{D}}$ is a torsion module as in the Theorem, then we have $ext^0(M)=hom_D(M,D)=0$ in any case.  \\

\noindent
{\bf REMARK 2.8}: Lanczos has been trying in vain to do for the $Bianchi$ operator what he did for the $Riemann$ operator, a useless  but possible 
"{\it shift by one step to the right} " and to do for the {\it Weyl} operator what he did for the {\it Riemann} operator. However, we shall discover that the dimension $n=4$, which is quite "{\it fine} " for the classical Killing sequence, is quite "{\it bad} " for the conformal Killing sequence, a result not known after one century because it cannot be understood without using the {\it Spencer $\delta$-comology} in the following commutative diagram which is explaining therefore what we shall call the "{\it Lanczos secret} ". Following ([33]) and the fact that the two central vertical $\delta$-sequences are exact, this diagram allows to construct the {\it Bianchi} operator ${\cal{D}}_2:F_1 \rightarrow F_2 $ as {\it generating} CC for the {\it Riemann} operator ${\cal{D}}_1: F_0=S_2T^* \rightarrow F_1$ defined by a similar diagram and thus only depends on the symbol $g_1$. We have the following commutative diagram allowing to define $F_1$, where all the rows are exact and the columns are exact but eventually the left one: \\

\[ \small { \begin{array}{rccccccccl}
   & 0 &  & 0  & &  0  &  &   &       \\
   & \downarrow &  &  \downarrow & & \downarrow & && &   \\
0 \rightarrow & g_3 & \rightarrow &S_3T^*\otimes T& \rightarrow &S_2T^ *\otimes F_0 &\rightarrow &F_1&\rightarrow & 0 \\
 & \downarrow  &  &  \downarrow & & \downarrow & & &   \\
0 \rightarrow & T^*\otimes g_2 & \rightarrow &T^*\otimes S_2T^*\otimes T& \rightarrow &T^*\otimes T^ *\otimes F_0 &\rightarrow & 0 &   \\
 & \downarrow  &  &  \downarrow & & \downarrow & & & &   \\
0 \rightarrow &{\wedge}^2 T^*\otimes g_1 & \rightarrow &\underline{{\wedge}^2T^*\otimes T^*\otimes T} & \rightarrow &{\wedge}^2T^*\otimes F_0 &\rightarrow & 0 && \\
& \downarrow  &  &  \downarrow & & \downarrow  & &  & &  \\
0 \rightarrow &{\wedge}^3 T^*\otimes T & = &{\wedge}^3T^*\otimes T& \rightarrow &0 & &  & &  \\
  & \downarrow  &  &  \downarrow & &  & & & &   \\
 
 & 0  &  & 0  & &    &  &   &   &   
\end{array} } \]   \\  

When $n=4$, we provide the respective fiber dimensions in the next diagram:  \\

 \[ \small { \begin{array}{rccccccccl}
   &  &  & 0  & &  0  &  &  &      \\
   &  &  &  \downarrow & & \downarrow & &  & &   \\
  & 0 &\rightarrow &80& \rightarrow & 100 &\rightarrow & 20&\rightarrow & 0  \\
 &   &  &  \downarrow & & \downarrow & & &   \\
   & 0 & \rightarrow & 160 & \rightarrow & 160 &\rightarrow & 0 & &    \\
 &  &  &  \downarrow & & \downarrow & & & &   \\
 0  \rightarrow &36 & \rightarrow & 96 &\rightarrow &60&\rightarrow & 0 &\\
& \downarrow &  &  \downarrow & & \downarrow & &  & &  \\
0 \rightarrow &16 & = &16 & \rightarrow &0 & & & &  \\
 & \downarrow  &  &  \downarrow & &  & & & &  \\
 & 0  &  & 0  & &    &  &   &   &   
\end{array} } \]   \\

\noindent
Using the {\it Spencer cohomology} at ${\wedge}^2T^*\otimes g_1$ and a {\it snake-type chase} ([35], p 174), the vector bundle $F_1=H^2(g_1)$ in 
this diagram or {\it Riemann candidate} in the language of Lanczos, is defined by the short exact sequence:  \\
\[ \begin{array}{rcccccl}
  0 \longrightarrow & F_1& \longrightarrow  & {\wedge}^2T^*\otimes g_1 & \stackrel{\delta}{\longrightarrow} & {\wedge}^3T^*\otimes T  &\longrightarrow 0  \\
  0 \longrightarrow &20 &\longrightarrow  & 36 & \stackrel{\delta}{\longrightarrow} & 16 & \longrightarrow 0  \
\end{array}   \]  
All the vertical down arrows are $\delta$-maps of Spencer and all the vertical columns are exact but the first, which may not be exact only at ${\wedge}^3T^*\otimes g_1$ with cohomology equal to $H^3(g_1)$ because we have:  \\
\[  g_1=\{{\xi}^k_i \in T^*\otimes T \mid {\omega}_{rj}{\xi}^r_i + {\omega}_{ir}{\xi}^r_j=0\} \simeq {\wedge}^2T^*\subset T^*\otimes T \,\, \stackrel{det(\omega)\neq 0}{\Longrightarrow}  \,\, g_2=0 \Rightarrow g_3=0 \Rightarrow  g_4=0  \]

Similarly, we obtain the following commutative diagram allowing to define $F_2$, with analogous comments:  \\

\[ \small { \begin{array}{rcccccccccl}
   & 0 &  & 0  & &  0  &  & 0  &  &     \\
   & \downarrow &  &  \downarrow & & \downarrow & & \downarrow & &  & \\
0 \rightarrow & g_4 & \rightarrow &S_4T^*\otimes T& \rightarrow &S_3T^ *\otimes F_0 &\rightarrow & T^*\otimes F_1&\rightarrow & F_2 & \rightarrow 0 \\
 & \downarrow  &  &  \downarrow & & \downarrow & & \parallel & &  \\
0 \rightarrow & T^*\otimes g_3 & \rightarrow &T^*\otimes S_3T^*\otimes T& \rightarrow &T^*\otimes S_2T^ *\otimes F_0 &\rightarrow & T^*\otimes F_1& \rightarrow & 0 &  \\
 & \downarrow  &  &  \downarrow & & \downarrow & & \downarrow & & &  \\
0 \rightarrow &{\wedge}^2 T^*\otimes g_2 & \rightarrow &{\wedge}^2T^*\otimes S_2T^*\otimes T& \rightarrow &{\wedge}^2T^*\otimes T^ *\otimes F_0 &\rightarrow &
 0&&& \\
& \downarrow  &  &  \downarrow & & \downarrow & &  & & & \\
0 \rightarrow &{\wedge}^3 T^*\otimes g_1 & \rightarrow &\underline{{\wedge}^3T^*\otimes T^*\otimes T}& \rightarrow &{\wedge}^3T^*\otimes F_0 &\rightarrow & 0 & &&  \\
  & \downarrow  &  &  \downarrow & & \downarrow & & & & &  \\
0 \rightarrow &{\wedge}^4 T^*\otimes T & = &{\wedge}^4T^*\otimes  T& \rightarrow &0 & & & &&  \\
 & \downarrow  &  &  \downarrow & &  & & & & & \\
 & 0  &  & 0  & &    &  &   &   &  & 
\end{array} } \]   \\  

When $n=4$, we provide the fiber dimensions in the next diagram:   \\

 \[ \small { \begin{array}{rcccccccccl}
   &  &  & 0  & &  0  &  & 0  &  &     \\
   &  &  &  \downarrow & & \downarrow & & \downarrow & &  & \\
  & 0 &\rightarrow &140& \rightarrow & 200 &\rightarrow & 80&\rightarrow & 20 & \rightarrow 0 \\
 &   &  &  \downarrow & & \downarrow & & \parallel & &  \\
   & 0 & \rightarrow & 320 & \rightarrow & 400 &\rightarrow & 80 & \rightarrow & 0 &  \\
 &  &  &  \downarrow & & \downarrow & & \downarrow & & &  \\
 &0 & \rightarrow &240 & \rightarrow & 240 &\rightarrow &0&&& \\
&  &  &  \downarrow & & \downarrow & &  & & & \\
0 \rightarrow &24 & \rightarrow & 64 & \rightarrow & 40 &\rightarrow & 0 & &&  \\
  & \downarrow  &  &  \downarrow & & \downarrow & & & & &  \\
0 \rightarrow &4& = &4& \rightarrow &0 & & & &&  \\
 & \downarrow  &  &  \downarrow & &  & & & & & \\
 & 0  &  & 0  & &    &  &   &   &  & 
\end{array} } \]   \\

\noindent
A {\it snake-type chase} similarly provides the identification  $F_2=H^3(g_1)$ while using again the {\it Spencer cohomology} at ${\wedge}^3T^*
\otimes g_1$. The vector bundle $F_2$ providing the Bianchi identities is thus defined by the exactness of the top row of the preceding diagram or, equivalently, using the left column, by the short exact sequence:  \\
\[  \begin{array}{rcccccl}
  0 \longrightarrow & F_2& \longrightarrow  & {\wedge}^3T^*\otimes g_1 & \stackrel{\delta}{\longrightarrow} & {\wedge}^4T^*\otimes T  &\longrightarrow 0  \\
  0 \longrightarrow &20 &\longrightarrow  &24 & \stackrel{\delta}{\longrightarrow} & 4 & \longrightarrow 0  \
\end{array}   \]  \\

When $n=4$, using the duality with respect to the volume form $dx^1\wedge dx^2\wedge dx^3 \wedge dx^4$ in order to change the indices, we obtain successively (care to the signs):  \\
\[  \begin{array}{rccccccccl}
   & B^i_{1,234} & - & B^i_{2,341} & + & B^i_{3,412} & - & B^i_{4,123} & = & 0  \\
   & B_{i1,1}  & - & B_{i2,2} & + & B_{i3,3} & - & B_{i4,4}&  =  &0   \\
i=4 \Rightarrow    &  B_{41,1}& - & B_{42,2}&+ &  B_{43,3} & = & 0 &   \\
  & L_{23,1}& + & L_{31,2}&  + &  L_{12,3} & = &  0
\end{array}  \]   \\  

\noindent
and finally exhibit the {\it Lanczos potential} $L\in {\wedge}^2T^*\otimes T^*$ as a $3$-tensor satisfying: \\

\noindent
\hspace*{3cm}  \fbox{$L_{ij,k} + L_{ji,k}=0, \,\,\,\,\, L_{ij,k} + L_{jk,i} + L_{ki,j}=0 $} \,\,\, $(24-4=20)$  \\  \\
but {\it this result is only valid in this specific situation} and does not provide any potential because {\it the adjoint sequence is going ... backwards} (!). \\

Starting with the (classical) {\it Killing} operator $K: T\rightarrow S_2T^*$ defined by $\xi \rightarrow {\cal{L}}(\xi)\omega $, we may obtain successively the following differential sequences for various useful dimensions:  \\

\noindent
\[ \begin{array}{rcccccccccccl}
n=2 \,\,\,\,\, \,\, & 2 & \underset 1{\stackrel{K}{\longrightarrow}} & 3 & \underset 2 {\stackrel{R}{\longrightarrow}}& 1 & \longrightarrow & 0 &  &  &  &  &\\
n=3 \,\,\,\,\,\,\,  & 3&   \underset 1{\stackrel{K}{\longrightarrow}} & 6& \underset 2 {\stackrel{R}{\longrightarrow}}& 6 & \underset 1{\stackrel{B}{\longrightarrow}} & 3 &\longrightarrow  
&  0 &  & & \\
n=4 \,\,\,\,\,\,\,  & 4 & \underset 1 {\stackrel{K}{\longrightarrow}} &10 & \underset 2{ \stackrel{R}{\longrightarrow}} & 20 & \underset 1 {\stackrel{B}{\longrightarrow}} & 20 & \underset 1 {\longrightarrow}& 6 & \longrightarrow & 0 &  \\
n=5 \,\,\,\,\,\,\,  & 5 & \underset 1 {\stackrel{K}{\longrightarrow}} &15 & \underset 2{ \stackrel{R}{\longrightarrow}} & 50 & \underset 1 {\stackrel{B}{\longrightarrow}} & 75 & \underset 1 {\longrightarrow}& 45 & \underset 1{\longrightarrow} &10 & \longrightarrow 0  
\end{array}   \] 
For example, we have the Euler-Poincar\'{e} characteristic:  $rk_D(M)=4-10+20-20+6=0$ when $n=4$ or $rk_D(M)=5 - 15 + 50 -75 + 45 - 10 = 0$ when $n=5$. Setting successively ${\cal{D}}= Killing, \,\, {\cal{D}}_1=Riemann, \, \, {\cal{D}}_2=Bianchi $ and so on, it follows from the previous study that each operator is parametrizing the following one. Applying {\it double duality} while introducing the respective adjoint operators, then $ad({\cal{D}}_2)$ is thus parametrizing  the {\it Beltrami} operator $ad(Riemann)=ad({\cal{D}}_1)$ with (canonical) potentials called {\it Lanczos} only when $n=4$ while $ad({\cal{D}}_1)$ is parametrizing the {\it Cauchy} operator $ad({\cal{D}})$ with (canonical) potentials called {\it Airy} when $n=2$, 
{\it Beltrami} when $n=3$ and so on ([28]). It must be finally noticed that $ad(Ricci)$ is also parametrizing the {\it Cauchy} operator while the 
{\it Einstein} operator is useless ([30]). \\ 

With more details, we now provide the explicit potentials and parametrizations of the adjoint sequence. In fact and up to our knowledge after more than twenty years of teaching continuum mechanics and elasticity ([31]), even the adjoint of the Killing operator is never presented within the differential duality and it is not completely evident that $ad(Killing)=Cauchy$, indepently of any EL constitutive relations, that is to say without refering to a Lagrangian. For this purpose, we recall that the standard {\it infinitesimal deformation tensor} is $\epsilon=\frac{1}{2}\Omega\in S_2T^*$, that is ${\epsilon}_{ij}= \frac{1}{2}{\Omega}_{ij}$. Accordingly, the {stress tensor density} is $\sigma \in {\wedge}^nT^*\otimes S_2T= ad(S_2T^*)$, a result leading most textbooks to conclude that the stress can be written as a symmetric matrix, contrary to its classical "{\it experimental} " presentation through the well known {\it Cauchy tetrahedral} that is {\it never} making any assumption on the symmetry of the underlying matrix. It is this result that pushed the brothers Cosserat to revisit the mathematical foundations of elasticity theory and to introduce a non-symmetrical stress {\it that could not have any relation with the above definition}. In actual practice, let us consider for simplicity the case $n=2$ with Euclidean metric. Then we have only $3$ (care) independent components of $\epsilon$, namely $({\epsilon}_{11}, {\epsilon}_{12}={\epsilon}_{21}, {\epsilon}_{22})$, that we should dualize with $({\sigma}^{11}, {\sigma}^{12}, {\sigma}^{22})$ and students know that the completion with ${\sigma}^{21}={\sigma}^{12}$, {\it automatically} done in EL ... and thus also GR, depends on a delicate proof involving equilibrium of torsors, a quite useful mechanical concept having no link with the previous procedure ([23],[24],[27]). \\
The "{\it only} " technical purpose is to arrive to a "{\it nice} " summation with factors $2$:  \\
\[ \begin{array}{rcl}
{\sigma}^{ij}{\epsilon}_{ij} & = &  {\sigma}^{11}{\epsilon}_{11} + {\sigma}^{12}{\epsilon}_{12} + {\sigma}^{21}{\epsilon}_{21} + {\sigma}^{22}{\epsilon}_{22} +  ...   \\
   &  =  &  {\sigma}^{11}{\epsilon}_{11} + \fbox{2} \,\, {\sigma}^{12}{\epsilon}_{12} + {\sigma}^{22}{\epsilon}_{22} +  ...   
\end{array}   \]
contrary to the "{\it pure} " duality sum:  \\
\[  {\sigma}^{11}{\epsilon}_{11} + {\sigma}^{12}{\epsilon}_{12} + {\sigma}^{22}{\epsilon}_{22} +  ...    \]
This apparently naive comment is in fact the deep reason for which the {\it Riemann} or {\it Beltrami} operators are self-adjoint $6\times 6$ operator matrices when $n=3$ while the {\it Einstein} operator is a self-adjoint $10\times 10$ operator matrix when $n=4$, contrary to the {\it Ricci} operator matrix. \\
Multiplying $\Omega \in F_0$ on the left by $\sigma \in ad(F_0)$ and integrating by parts the $n$-form, we get:  \\
\[  \begin{array}{rcl}
\frac{1}{2} {\sigma}^{ij}{\Omega}_{ij} & = & \frac{1}{2} {\sigma}^{ij}({\omega}_{rj}{\partial}_i{\xi}^r + {\omega}_{ir} {\partial}_j{\xi}^r +{\xi}^r {\partial}_r{\omega}_{ij})  \\
  &  =  & - ({\partial}_i{\sigma}^{ri} + \frac{1}{2} ({\partial}_i{\omega}_{rj} + {\partial}_j{\omega}_{ir} - 
  {\partial}_r{\omega}_{ij}){\sigma}^{ij}){\xi}^r + ... \\
    &  =  &  - ({\partial }_i{\sigma}^{ri} + {\gamma}^r_{ij} {\sigma}^{ij}){\xi}_r  + ... \\
    &  =  &  - {\nabla}_i{\sigma}^{ri}{\xi}_r + ... 
\end{array}   \]
where the covariant derivative takes into account  the tensorial density nature of $\sigma$. This operator is reducing locally to the classical {\it Cauchy} equations $d_r{\sigma}^{ir}=f^i$ by introducing the force density $f\in {\wedge}^nT^*\otimes  T^*=ad(T)$ in the right member and lowering the index by means of $\omega$ as usual. \\

Exactly the same procedure can be applied to EM while starting with the {\it field equations} (first set of Maxwell equations) $dF=0$ where 
$F\in {\wedge}^2T^*$ and $d:{\wedge}^2T^* \rightarrow {\wedge}^3T^*$ is the standard exterior derivative. Using the local exactness of the 
Poincar\'{e} sequence, we may find a parametrization $dA=F$ with $A\in T^*$. The {\it induction equations} (second set of Maxwell equations) are described, indpendently of any Lagrangian or EM Minkowski constitutive relations, by (See [31] for detais):  \\
\[  ad(d): {\wedge}^4T^* \otimes {\wedge}^2T \rightarrow {\wedge}^4T^*\otimes T:{\cal{F}}\rightarrow {\cal{J}}\Leftrightarrow {\partial}_r{\cal{F}}^{ir}={\cal{J}}^i  \]
Both sets of Maxwell equations are invariant by any local diffeomorphism and the conformal group of space-time is {\it only} the biggest group of invariance of the Minkowski constitutive laws in vacuum ([18],[19,[26],[27]). Of course, this result is showing that, contrary to the existence of the well known EL/EM couplings (Piezoelectricity, Photoelasticity, streaming birefringence, ...) where $\Omega$ {\it and} $F$ should appear equally in Lagrangians and constitutive relations, there is no room for the EM field $F$ in the Janet sequence and no room for the EL field $\Omega$ in the Spencer sequence. It follows that there cannot be any relation existing between the EM field $F$ and the Riemann tensor $R$, contrary to what Lanczos was believing.  \\

Coming back to the initial resolution of Killing vector fields and its adjoint, we must push by one step to the right in order to study the {\it Riemann} operator and its adjoint, the {\it Beltrami} operator. Of course, from its definition, the $Riemann$ operator is parametrized ... by the {\it Killing}  operator. As the sequence can only be constructed for metrics with constant curvature, thus for flat metrics, we assume that we deal with Euclidean metric for EL when $n\geq 2$ or with the Minkowski metric for GR when $n=4$ as a way to simplify the formulas (See [33] otherwise). In order to correct the formula $(III.6)$ of ([11]) with no mathematical meaning, we may linearize ($\omega$, $\gamma$, $\rho$) and obtain successively:  \\
\[   2 {\Gamma}^k_{ij}= {\omega}^{kr}(d_i{\Omega}_{rj} + d_j{\Omega}_{ir} - d_r{\Omega}_{ij}) , \,\,\,\,\,
R^k_{l,ij} = d_i{\Gamma}^k_{lj} - d_j {\Gamma}^k_{li} \] 
\[\Rightarrow  \,\,  2\,\, R_{kl,ij}=(d_{li}{\Omega}_{kj} - d_{lj}{\Omega}_{ki}) - (d_{ki}{\Omega}_{lj} - d_{kj}{\Omega}_{li})   \]
Then we check, like in the Introduction, that such an $R$ is a section of the {Riemann} candidate $F_1$. In order to understand the underlying confusion, we let the reader prove as an exercise that the number of CC of the second order system $d_{ij}{\xi}^k={\Gamma}^k_{ij}$, namely $d_i{\Gamma}^k_{lj} -d_j{\Gamma}^k_{li}=0, \forall 1\leq i,j,r\leq n $, is equal to $n^2(n^2-1)/12$, that is ... $120$ when $n=4$ and {\it not}  $20$.  \\
We shall use the dual notations:  \\
\[  \begin{array}{ccccccc}
\xi & \rightarrow & \Omega & \rightarrow & R & \rightarrow & B  \\
f  & \leftarrow &\sigma & \leftarrow & \alpha & \leftarrow & \lambda
\end{array}   \]
Multiplying $R \in F_1$ on the left by $\alpha \in ad(F_1)$ and integrating by parts the $n$-form, we get:  \\
\[  \begin{array}{rcl}
 2 \,\alpha R & = & {\alpha}^{kl,ij}(d_{li}{\Omega}_{kj} - d_{ki}{\Omega}_{lj} + d_{kj}{\Omega}_{li} - d_{lj}{\Omega}_{ki}) \\
  &  = & (d_{li} {\alpha}^{kl,ij}){\Omega}_{kj} -(d_{ki}{\alpha}^{kl,ij}){\Omega}_{lj} + 
(d_{kj}{\alpha}^{kl,ij}) {\Omega}_{li} - (d_{lj}{\alpha}^{kl,ij}){\Omega}_{ki}  + ...  \\
 &  =  &  (d_{rs}{\alpha}^{ir,sj} - d_{rs}{\alpha}^{ri,sj} + d_{rs}{\alpha}^{ri,js} - d_{rs}{\alpha}^{ir,js} ){\Omega}_{ij} + ...  \\
  &  =  &  4\,\,(d_{rs}{\alpha}^{ir,sj}){\Omega}_{ij} +  ...
\end{array}  \]
and the striking parametrization $d_{rs}{\alpha}^{ir,sj} = {\sigma}^{ij}={\sigma}^{ji} \Rightarrow d_j{\sigma}^{ij}=0 $. When $n=2$, setting 
$ {\alpha}^{12,12}= - \phi$, we get the Airy parametrization ${\sigma}^{11}= d_{22}\phi, {\sigma}^{12}={\sigma}^{21}= - d_{12}\phi, {\sigma}^{22}= d_{11}\phi$ ([23]).  \\

Finally, we shall construct the {\it Lanczos} operator as the adjoint of the {\it Bianchi} operator and will thus justify the formula $(III.5)$ in ([11]) when $n=4$:  \\
\[   (d_rR_{kl,ij} + d_iR_{kljr} + d_jR_{klri}=B_{kl,ijr}) \in {\wedge}^3T^*\otimes g_1\simeq {\wedge}^3T^*\otimes {\wedge}^2T^*  \]
Multiplying $B \in F_2$ on the left by $\lambda \in ad(F_2)$ and integrating by parts the $n$-form, we get:  \\
\[   \begin{array}{rcl}
   \lambda B & = & {\lambda}^{kl,ijr}(d_rR_{kl,ij} + d_iR_{kljr} + d_jR_{klri})   \\
                  &  = & - ((d_r{\lambda}^{kl,ijr}) R_{kl,ij} + (d_i{\lambda}^{kl,ijr})R_{kl,jr} +(d_j{\lambda}^{kl,ijr})R_{kl,ri})  +  ... \\
                  &  =  &  -3\,\, (d_r{\lambda}^{kl,ijr}) R_{kl,ij} +  ...
              \end{array}   \]
and the striking parametrization  $ d_r{\lambda}^{kl,ijr}={\alpha}^{kl,ij}\Rightarrow d_{li}{\alpha}^{kl,ij}=d_{rli}{\lambda}^{kl,ijr}=0$ up to sign. This is the main result explainging the confusion done by Lanczos between Hodge duality and differential duality, that is between the {\it Killing} sequence and its adjoint which is going ... {\it backwards }!. \\
We have thus confirmed the fact that the differential module defined by the {\it Killing} operator ${\cal{D}}$ has vanishing first and second differential extension modules while its adjoint differential module defined by $ad(Killing)=Cauchy$ is a reflexive differential module. The reader will have noticed that not a single of the previous computations could be even imagined without these new tools which are at the same time explicit and intrinsic.  \\

\vspace{15mm}

\noindent
{\bf 3) WEYL/LANCZOS POTENTIAL}  \\

Starting now afresh with the {\it conformal Killing} operator $CK$ such that ${\cal{L}}(\xi)\omega=A(x)\omega $ or, equivalently, introducing the {\it metric density} ${\hat{\omega}}_{ij}={\omega}_{ij} {\mid det(\omega)\mid}^{ - \frac{1}{n}}$, we have a new operator $CK: T \rightarrow \{\Omega \in S_2T^*\mid tr(\Omega)={\omega}^{ij}{\Omega}_{ij}=0\}$ defined by $\xi \rightarrow {\cal{L}}(\xi)\hat{\omega}$ and we obtain successively the following differential sequences for various useful dimensions:  \\

\[ \begin{array}{rcccccccccccl}
n=3 \,\,\,\,\,\,\,  & 3 & \underset 1{\stackrel{CK}{\longrightarrow}} & 5 & \underset 3 {\stackrel{?}{\longrightarrow}}& 5 & \underset 1 {\stackrel{?}{\longrightarrow}} & 3 & \longrightarrow & 0 & &  &  \\
n=4 \,\,\,\,\, \,\, & 4&   \underset 1{\stackrel{CK}{\longrightarrow}} & 9 & \underset 2 {\stackrel{W}{\longrightarrow}} & 10 & \underset 2{\stackrel{CB}{\longrightarrow}} & 9 &\underset 1{\longrightarrow } &  4 &\longrightarrow  &0  & \\
n=5 \,\,\,\,\, \,\, & 5 & \underset 1 {\stackrel{CK}{\longrightarrow}} &14 & \underset 2{ \stackrel{W}{\longrightarrow}} & 35 & \underset 1 
{\stackrel{CB}{\longrightarrow}} & 35 & \underset 2 {\longrightarrow}& 14 & \underset 1 {\longrightarrow} &5 &\longrightarrow 0
\end{array}   \]  
 
For example, we have the Euler-Poincar\'{e} characteristic:  $rk_D(M)=5-14+35-35+14-5=0$. These results have been confirmed by computer algebra in ([29]). They prove that the analogue of the {\it Weyl} operator is of order $3$ when $n=3$ but becomes of order $2$ when $n\geq 4$ and that the analogue of the 
{\it Bianchi} operator is now of order $1$ when $n=3$, of order $2$ when $n=4$ but becomes again of order $1$ when $n\geq 5$. \\  
Proceeding exactly as before in order to define ${\hat{F}}_2$, we obtain the following commutative diagram where all the rows are exact and the columns are exact but eventually the left one: \\

\noindent
           \[ \small { \begin{array}{rcccccccccl}
   & 0 &  & 0  & &  0  &  & 0  &  &     \\
   & \downarrow &  &  \downarrow & & \downarrow & & \downarrow & &  & \\
0 \rightarrow & {\hat{g}}_5 & \rightarrow &S_5T^*\otimes T& \rightarrow &S_4T^ *\otimes {\hat{F}}_0 &\rightarrow & S_2T^*\otimes 
{\hat{F}}_1&\rightarrow & {\hat{F}}_2 & \rightarrow 0 \\
 & \downarrow  &  &  \downarrow & & \downarrow & & \downarrow & &  \\
0 \rightarrow & T^*\otimes {\hat{g}}_4 & \rightarrow &T^*\otimes S_4T^*\otimes T& \rightarrow &T^*\otimes S_3T^ *\otimes {\hat{F}}_0 &\rightarrow & T^*\otimes T^*\otimes {\hat{F}}_1& \rightarrow & 0 &  \\
 & \downarrow  &  &  \downarrow & & \downarrow & & \downarrow & & &  \\
0 \rightarrow &{\wedge}^2 T^*\otimes {\hat{g}}_3 & \rightarrow & {\wedge}^2T^*\otimes S_3T^*\otimes T& \rightarrow &{\wedge}^2T^*\otimes S_2T^ *\otimes {\hat{F}}_0 &\rightarrow &{\wedge}^2T^*\otimes {\hat{F}}_1 & \rightarrow & 0 & \\
& \downarrow  &  &  \downarrow & & \downarrow & & \downarrow & & & \\
0 \rightarrow &{\wedge}^3 T^*\otimes {\hat{g}}_2 & \rightarrow &\underline{{\wedge}^3T^*\otimes S_2T^*\otimes T}& \rightarrow &{\wedge}^3T^*\otimes T^*\otimes {\hat{F}}_0 &\rightarrow & 0 & &&  \\
  & \downarrow  &  &  \downarrow & & \downarrow & & & & &  \\
0 \rightarrow &{\wedge}^4 T^*\otimes {\hat{g}}_1 & \rightarrow &{\wedge}^4T^*\otimes T^*\otimes  T& \rightarrow & {\wedge}^4T^*
\otimes {\hat{F}}_0 &\rightarrow  & 0
& &&  \\
 & \downarrow  &  &  \downarrow & & \downarrow   & & & & & \\
 & 0  &  & 0  & &   0  &  &   &   &  & 
\end{array} } \]   \\  

When $n=4$, we provide the fiber dimensions in the next diagram:   \\

\noindent
           \[ \small { \begin{array}{rcccccccccl}
   &  &  & 0  & &  0  &  & 0  &  &     \\
   &  &  &  \downarrow & & \downarrow & & \downarrow & &  & \\
 & 0 & \rightarrow &224& \rightarrow &315 &\rightarrow & 100 &\rightarrow & 9 & \rightarrow 0 \\
 &  &  &  \downarrow & & \downarrow & & \downarrow & &  \\
 & 0 & \rightarrow & 560 & \rightarrow & 720 &\rightarrow & 160 & \rightarrow & 0 &  \\
 &  &  &  \downarrow & & \downarrow & & \downarrow & & &  \\
 & 0  & \rightarrow &  480 & \rightarrow & 540 &\rightarrow & 60 & \rightarrow & 0 & \\
&  &  &  \downarrow & & \downarrow & & \downarrow & & & \\
0 \rightarrow & 16 & \rightarrow & 160 & \rightarrow &144 &\rightarrow & 0 & &&  \\
  & \downarrow  &  &  \downarrow & & \downarrow & & & & &  \\
0 \rightarrow &7 & \rightarrow & 16 & \rightarrow & 9 &\rightarrow  & 0
& &&  \\
 & \downarrow  &  &  \downarrow & & \downarrow   & & & & & \\
 & 0  &  & 0  & &   0  &  &   &   &  & 
\end{array} } \]   \\  

It is much more difficult to prove that the last map $\delta: {\wedge}^3T^*\otimes {\hat{g}}_2 \rightarrow {\wedge}^4T^*\otimes T: 16 \rightarrow 7$ is an epimorphism. We let the reader manage as an exercise through an up and down delicate circular chase in order to convince him that {\it no classical result could provide such a result} which is nevertheless an {\it obligatory step} for finding the desired $dim({\hat{F}}_2)= 16- 7=9$ ({\it Hint}: Prove first that the $Weyl$ operator ${\hat{\cal{D}}}_1$ has no first order CC when $n=4$, then prove that each element of ${\wedge}^2T^*\otimes S_2T^*\otimes {\hat{F}}_0$ is the sum of an element in $\delta (T^*\otimes S_3T^*\otimes {\hat{F}}_0)$ and an element coming from ${\wedge}^2T^*\otimes S_3T^*\otimes T$). Of course, in view of the dimensions of the matrices involved (up to $540 \times 720$), we wish good luck to anybody trying to use computer algebra and refer to the computations done in ([29]) that have been done while knowing "{\it a priori} " the dimensions that should be found.   \\

We finally notice that the change of the successive orders is totally unusual as in ([33]) and refer to ([32]) for more details on the computer algebra methods. In particular, when $n=4$, the conformal analogue of the {\it Bianchi} operator is now of order $2$, a result explaining why Lanczos and followers could not succeed adapting the Lanczos tensor potential $L$ for the {\it Weyl} operator, even if it was already known for the {\it Riemann} operator when $n=4$. In particular, thanks to Theorems 2.1 and Corollary 2.4, we have thus solved the {\it Riemann-Lanczos} and {\it Weyl-Lanczos} parametrization problems in arbitrary dimension while providing explicit computations.\\

\newpage

\noindent
{\bf 4) RIEMANN VERSUS WEYL}  \\

\noindent
Using the splitting of the Riemann tensor between the Ricci or Einstein tensor and the Weyl tensor in the second column while taking into account the fact that the extension modules are torsion modules, then each component of the Weyl tensor is differentially dependent on the Ricci tensor and we recall the commutative and exact diagram first provided in ([30]) when $n=4$:  \\
\[  \begin{array}{cccccccl}
  &  &  0  &  &  0  &  &  0  &     \\
  & &  \downarrow & &  \downarrow  &  &  \downarrow  &  \\
  0 &  & 10 & \longrightarrow & 16 & \rightarrow & 6 &  \rightarrow 0  \\
  \downarrow &  &  \downarrow \uparrow&  & \downarrow &  & \parallel &   \\
  10  &  \stackrel{Riemann}{\longrightarrow} & 20 & \stackrel{Bianchi}{\longrightarrow} & 20 & \rightarrow & 6 & \rightarrow 0  \\
\parallel  &  &  \downarrow\uparrow &  & \downarrow &  & \downarrow &   \\
10 & \stackrel{Einstein}{\longrightarrow}&  10 &  \stackrel{div}{\longrightarrow} & 4  & \rightarrow & 0 &  \\
\downarrow &  & \downarrow & & \downarrow & & &  \\
0 & & 0 & & 0 & & & 
\end{array}   \]
It follows that the $10$ components of the Weyl tensor must satisfy a first order linear system with $16$ equations, having $6$ generating first order CC. The differential rank of the corresponding operator is thus equal to the differential rank of its image that is $16-6=10$ and such an operator defines therefore a torsion module because the differential rank of its kernel is $10-10=0$. Equivalently, we have to look {\it separately} for each component of the Weyl tensor in order to obtain the Lichnerowicz {\it wave equations} (as they are called in France !) ([[30]). The situation is similar to that of the Cauchy-Riemann equations obtained when $n=2$ by considering the {\it conformal Killing} operator. Indeed, any complex transformation $y=f(x)$ must be solution of the (linear) first order system $y^2_2 - y^1_1=0, y^1_2 + y^2_1=0$ of finite Lie equations though we obtain $y^1_{11} +y^1_{22}=0, y^2_{11}+y^2_{22}=0$, that is $y^1$ and $y^2$ are {\it separately} killed by the second order {\it Laplace} operator $\Delta =d_{11}+d_{22}$. We obtain the following striking technical lemma explaining the so-called {\it gauging} procedure of the Lanczos potential. \\

\noindent
{\bf LEMMA 4.1}: When $n=4$, the central vertical arrow $20 \rightarrow 4$ is just described by the contraction formula:   \\
\[ L=L_{ij,k}\,\,\, \rightarrow \,\,\,  L_i={\omega}^{jk}L_{ij,k} \] \\
\noindent
{\it Proof}: Let us write down the $Bianchi$ operator in the form:   \\
\[  {\nabla}_r{R^k}_{l,ij} + {\nabla}_i{R^k}_{l,jr} + {\nabla}_j{R^k}_{l,ri}= {B^k}_{l,ijr}   \]
Contracting with $k=j=s$, we obtain:  \\
\[   {\nabla}_r{R^s}_{l,is} + {\nabla}_i{R^s}_{l,sr} + {\nabla}_s{R^s}_{l,ri} = {B^s}_{l,isr}   \]
Setting as usual ${R^s}_{l,sr}=R_{lr}=R_{rl}$ with ${\omega}^{ij}R_{ij}=R$ and contracting with ${\omega}^{li}$, we finally get :   \\
\[  2 {\nabla}_s {R^s}_r - {\nabla}_r R= {\omega}^{ij}{B^s}_{i,jrs}   \]
as {\it the} way to use a contraction in order to exhibit Einstein equations.  \\
With $n=4$, let us write down all the terms, using the Euclidean metric for simplicity instead of the Minkowski metric, recalling that only this later choice allows to find out both the Poincar\'{e} group {\it and} the differential sequence with successive operators $K,R,B$ according to ([33]):  \\
\[  {B^s}_{1,1rs} + {B^s}_{2,2rs} + {B^s}_{3,3rs} + {B^s}_{4,4,rs} = C_r    \]
that is to say with all the terms:  \\
\[ { \left\{ \begin{array}{c}
\,\,\, {B^1}_{1,1r1} + {B^2}_{1,1r2} + {B^3}_{1,1,r3} + {B^4}_{1,1r4} \\
+{B^1}_{2,2r1} + {B^2}_{2,2r2} + {B^3}_{2,2,r3} + {B^4}_{2,2r4}\\
+{B^1}_{3,3r1} + {B^2}_{3,3r2} + {B^3}_{3,3,r3} + {B^4}_{3,3r4} \\
+{B^1}_{4,4r1} + {B^2}_{4,4r2} + {B^3}_{4,4r3} + {B^4}_{4,4r4} 
\end{array}\right. }  = C_r  \]
where, {\it in any case}, we have ${B^1}_{1,1r1}={B^2}_{2,r2}={B^3}_{3,r3}={B^4}_{4,r4}=0 $.  \\
If we set $r=1$, the first line disappears because of the $3$-form ${\wedge}^3T^*$ and we are left with:  \\
\[  {  \left\{  \begin{array}{c}
\,\,\, {B^3}_{2,213} + {B^4}_{2,214}  \\
+{B^2}_{3,312} + {B^4}_{3,3,14}  \\
+{B^2}_{4,412} + {B^3}_{4,413}
\end{array}  \right.  }  = C_1  \]
Using Hodge duality, we get with new indices:  \\
\[  {  \left\{  \begin{array}{c}
- {B^3}_{2,4} + {B^4}_{2,3}  \\
+ {B^2}_{3,4}  - {B^4}_{3,2}  \\
- {B^2}_{4,3} + {B^3}_{4,2}
\end{array}  \right.  } = C_1  \]
arriving finally to the formula:  \\
\[     2({B^3}_{4,2} + {B^4}_{2,3} + {B^2}_{3,4})=C_1  \]
that is {\it exactly} twice the trace of the Lanczos tensor, namely:  \\
\[   {{L_1}^2}_2 +  {{L_1}^3}_3 +{{L_1}^4}_4 = {{L_1}^r}_r   \]
This result explains why the Lanczos tensor $L_{ij,k}= - L_{ji,k}$ with $24$ components is first reduced to $20$ components through the condition 
$L_{ij,k} + L_{jk,i} + L_{k,ij}=0$ and finally to $16$ components as in the diagram through the kernel of the above trace condition. It is thus {\it impossible} to understand this result even for $n=4$ without the Spencer $\delta$-cohomology and {\it absolutely impossible} to generalize this result in arbitrary dimension without the combination of the $\delta$-cohomology {\it and} double duality in differential homological algebra.  \\
\noindent
\hspace*{12cm}  Q.E.D.   \\

Using the previous definition $ad(E)={\wedge}^nT^*\otimes E^*$, such a result explains the confusion done by Lanczos and followers between the 
{\it Riemann} candidate $F_2$ or the {\it Weyl} candidate ${\hat{F}}_2$ and their respective formal adjoint vector bundles having of course the same fiber dimension but different transition rules under changes of local coordinates.  \\  \\

\vspace{2cm}

\noindent
{\bf CONCLUSION}  \\

When there is a competition between mathematics (differential homological algebra, double duality) and physics (Einstein equations, Maxwell equations), coming from their mixing up, sooner or later mathematics is always winning. This has been typically the situation met with the Lanczos potential theory where the motivating idea was quite clever but the final achievement has been contradictory with group theory through the only introduction of the Riemann tensor and Einstein equations, but without any reference to the conformal group of space-time and to the Weyl tensor that does not seem to have been known by Lanczos, even as late as in 1967. As we explained in the Introduction, the reader must nevertheless not forget that that it was not possible to discover any solution of the parametrization problem by potentials through differential double duality before $1990/1995$, that is too late for the many people already engaged in this type of research. We have clarified the situation with the powerful mathematical tools existing today and hope that computer algebra will take profit of this fact in the future.  \\  \\

\newpage

\noindent
{\bf REFERENCES}  \\

\noindent
[1] E. CARTAN: Les Syst\`{e}mes Diff\'{e}rentiels Ext\'{e}rieurs et Leurs Applications G\'{e}om\'{e}triques, Hermann, Paris,1945.  \\
\noindent 
[2]  P. DOLAN, A. GERBER: Janet-Riquier Theory and the Riemann-Lanczos Problem in $2$ and $3$ Dimensions, 2002, arXiv:gr-gq/0212055.  \\
\noindent
[3] P. DOLAN, A. GERBER: The Riemann-Lanczos Problem as an Exterior Differential Systm with Examples in $4$ and $5$ Dimensions, 2002, arXiv.gr-qc/0212054.  \\
\noindent
[4] S.B. EDGAR: On Effective Constraints for the Riemann-Lanczos Systems of Equations. J. Math. Phys., 44 (2003) 5375-5385. \\ 
 http://arxiv.org/abs/gr-qc/0302014   \\
\noindent
[5] S.B. EDGAR, A. H\"{O}GLUND: The Lanczos potential for Weyl-Candidate Tensors Exists only in Four Dimension, General Relativity and Gravitation, 32, 12 (2000) 2307. \\
http://rspa.royalsocietypublishing.org/content/royprsa/453/1959/835.full.pdf   \\ 
 \noindent
[6] S.B. EDGAR, J.M.M. SENOVILLA: A Local Potential for the Weyl tensor in all dimensions, Classical and Quantum Gravity, 21 (2004) L133.\\
 http://arxiv.org/abs/gr-qc/0408071     \\
\noindent
[7] L.P. EISENHART: Riemannian Geometry, Princeton University Press, 1926.  \\
\noindent
[8] M. JANET: Sur les Syst\`{e}mes aux D\'{e}riv\'{e}es Partielles, Journal de Math., 8(3) (1920) 65-151.  \\
\noindent
[9]  C.LANCZOS: A Remarkable Property of the Riemann-Christoffel Tensor in Four Dimensions, Annals of Math., 39 (1938) 842-850. \\
\noindent
[10] C. LANCZOS: Lagrange Multiplier and Riemannian Spaces, Reviews of Modern Physics, 21 (1949) 497-502.  \\ 
\noindent
[11] C. LANCZOS: The Variation Principles of Mechanics, Dover, New York, 4th edition, 1949.  \\
\noindent
[12] C. LANCZOS: The Splitting of the Riemann Tensor, Rev. Mod. Phys. 34, 1962, 379-389.  \\
\noindent
[13] C. LANCZOS: Le Tenseur de Riemann \`{a} Quatre Dimensions, Annales Scientifiques de l'Universit\'{e} de Clermont-Ferrand 2, Math.,Tome 8, 2 (1962) 167-170.  \\
  http://numdam.org/item/ASCFM\underline{\hspace{1.5mm}}1962\underline{\hspace{3mm}}8\underline{\hspace{1.5mm}}2\underline{\hspace{1.5mm}}167\underline{\hspace{1.5mm}}0         \\
\noindent
[14] U. OBERST: Multidimensional Constant Linear Systems, Acta Appl. Math., 20 (1990) 1-175.   \\
\noindent
[15] P. O'DONNELL, H. PYE: A Brief Historical Review of the Important Developments in Lanczos Potential Theory, EJTP, 24 (2010) 327-350.  \\
\noindent
[16] J.-F. POMMARET: Systems of Partial Differential Equations and Lie Pseudogroups, Gordon and Breach, New York, 1978 
(Russian translation by MIR, Moscow, 1983) \\
\noindent
[17]  J.-F. POMMARET: Differential Galois Theory, Gordon and Breach, New York, 1983.  \\
\noindent
[18]  J.-F. POMMARET: Lie Pseudogroups and Mechanics, Gordon and Breach, New York, 1988.\\
\noindent
[19] J.-F. POMMARET: Partial Differential Equations and Group Theory, New Perspectives for Applications, Mathematics and its Applications 293, Kluwer, 1994.\\
http://dx.doi.org/10.1007/978-94-017-2539-2   \\
\noindent
[20] J.-F. POMMARET: Dualit\'{e} Diff\'{e}rentielle et Applications, C. R. Acad. Sci. Paris, 320, S\'{e}rie I (1995) 1225-1230.  \\
\noindent
[21] J.-F. POMMARET: Partial Differential Control Theory, Kluwer, 2001, 957 pp.\\
\noindent
[22] J.-F. POMMARET: Algebraic Analysis of Control Systems Defined by Partial Differential Equations, in Advanced Topics in Control Systems Theory, Lecture Notes in Control and Information Sciences 311, Chapter 5, Springer, 2005, 155-223.\\
\noindent
[23] J.-F. POMMARET: Parametrization of Cosserat Equations, Acta Mechanica, 215 (2010) 43-55.\\
\noindent
[24] J.-F. POMMARET: Spencer Operator and Applications: From Continuum Mechanics to Mathematical Physics, in "Continuum Mechanics-Progress in Fundamentals and Engineering Applications", Dr. Yong Gan (Ed.), ISBN: 978-953-51-0447--6, InTech, 2012, Chapter 1: \\
http://www.intechopen.com/books/continuum-mechanics-progress-in-fundamentals-\\and-engineering-applications   \\
\noindent
[25] J.-F. POMMARET: The Mathematical Foundations of General Relativity Revisited, Journal of Modern Physics, 4 (2013) 223-239.\\
http://dx.doi.org/10.4236/jmp.2013.48A022   \\
\noindent
[26] J.-F. POMMARET: The Mathematical Foundations of Gauge Theory Revisited, Journal of Modern Physics, 5 (2014) 157-170.  \\
http://dx.doi.org/10.4236/jmp.2014.55026    \\
\noindent
[27] J.-F. POMMARET: From Thermodynamics to Gauge Theory: the Virial Theorem Revisited, in " Gauge Theories and Differential geometry ", NOVA Science Publishers, 2015, Chapter 1, 1-44.  \\
http://arxiv.org/abs/1504.04118  \\
\noindent
[28] J.-F. POMMARET: Airy, Beltrami, Maxwell, Einstein and Lanczos Potentials Revisited, Journal of Modern Physics, 7 (2016) 699-728.  \\
http://dx.doi.org/10.4236/jmp.2016.77068   \\
\noindent
[29] J.-F. POMMARET: Deformation Theory of Algebraic and Geometric Structures, Lambert Academic Publisher, (LAP), Saarbrucken, Germany, 2016.  \\
http://arxiv.org/abs/1207.1964  \\
\noindent
[30] J.-F. POMMARET: Why Gravitational Waves Cannot Exist, Journal of Modern Physics, 8,13 (2017) 2122-2158.  \\
http://dx.doi.org/10.4236/jmp.2017.813130   \\
\noindent
[31] J.-F. POMMARET: From Elasticity to Electromagnetism: Beyond the Mirror, \\
http://arxiv.org/abs/1802.02430  \\
\noindent
[32] J.-F. POMMARET: New Mathematical Methods for Physics, NOVA Science Publisher, New York, 2018.  \\
\noindent
[33] J.-F. POMMARET: Minkowski, Schwarzschild and Kerr Metrics Revisited, Journal of Modern Physics, 9 (2018) 1970-2007.  \\
https://doi.org/10.4236/jmp.2018.910125   \\
arXiv:1805.11958v2   \\
\noindent
[34] J.-F. POMMARET: Homological Solution of the Riemann-Lanczos and Weyl-Lanczos Problems in Arbitrary Dimension, 2018.  \\
http://arxiv.org/abs/1803.09610  \\
\noindent
[35] J. J. ROTMAN: An Introduction to Homological Algebra, Academic Press, 1979.\\
\noindent
[36] J.-P. SCHNEIDERS: An Introduction to D-Modules, Bull. Soc. Roy. Sci. Li\`{e}ge, 63 (1994) 223-295.  \\
\noindent
[37] D.C. SPENCER: Overdetermined Systems of Partial Differential Equations, Bull. Amer. Math. Soc.,75 (1965) 1-114.  \\
\noindent
[38] E. VESSIOT: Sur la Th\'{e}orie des Groupes Infinis, Ann. Ec. Normale Sup., 20 (1903) 411-451 (Can be obtained from  http://numdam.org).  \\

\end{document}